IAC-24-D2,8,x85968

# OPTIMISATION OF CATEGORICAL CHOICES IN EXPLORATION MISSION CONCEPTS OF OPERATIONS USING COLUMN GENERATION METHOD

**Nicholas Gollins [a*], Masafumi Isaji [a], and Koki Ho [a]**

[a] *Georgia Institute of Technology, Atlanta, Georgia, United States, 30332,* ngollins3@gatech.edu
\* Corresponding Author

**Abstract**
Space missions, particularly complex, large-scale exploration campaigns, can often involve a large number of discrete decisions or events in their concepts of operations. Whilst a variety of methods exist for the optimisation of continuous variables in mission design, the inherent presence of discrete events in mission ConOps disrupts the possibility of using methods that are dependent on having well-defined, continuous mathematical expressions to define the systems and instead creates a categorical mixed-integer problem. Typically, mission architects will circumvent this problem by solving the system optimisation for every permutation of the categorical decisions if practical, or use metaheuristic solvers if not. However, this can be prohibitively expensive in terms of computation time. Alternatively, categorical decisions in optimisation problems can be expressed using binary variables that indicate if the decision was taken or not. If implemented naively, commercially available mixed integer linear optimisation solvers are still slow to solve such a problem, in some cases not performing much better than combinatorially testing every permutation of the ConOps. Problems of this class can be solved more efficiently using "column generation" methods. Here, smaller, simpler restricted problems are created by removing significant numbers of variables. The restricted problem is solved, and the unused variables are priced by examining the dual linear program in order to test which, if any, could improve the objective of the restricted problem if they were to be added. Column generation methods are problem-specific, and so there is no guaranteed solution to these categorical problems. As such, the following paper proposes guidelines for defining restricted problems representing space exploration mission concepts of operations featuring common categories of decisions. First, the column generation process is described and then applied to two case studies. Firstly, it is applied to the NASA Marshall Advanced Concepts Office (ACO) ConOps for a crewed Mars mission, in which the design, assembly, and staging of the trans-Martian spacecraft are modelled using discrete decisions. Secondly, the process is applied to the payload delivery scheduling of translunar logistics in the context of an extended Artemis surface exploration campaign model.

**Nomenclature**
$A$ = matrix of constraint coefficients
$b$ = vector of constraint bounds
$B$ = scheduling binary variable (Case Study 2)
$c$ = vector of cost coefficients (Background and Theory), commodity time (Case Study 2)
$\mathcal{C}$ = set of co-manifested payloads
$d$ = demand matrix (Case Study 2)
$\mathcal{D}$ = domain (allowed set of arcs of a vehicle) (Case Study 2)
$f$ = primal objective function
$g$ = constraint function
$h$ = dual objective function
$I_{sp}$ = specific impulse
$k$ = crew consumable consumption rate (Case Study 2)
$\mathcal{L}$ = Lagrangian function
$m$ = mass
$N$ = maximum number of manoeuvres (Case Study 1), dimensions of index (Case Study 2)
$o$ = parking orbit selection binary variable
$\mathcal{O}$ = set of available parking orbits
$\mathcal{P}$ = set of "soft" precursor payloads
$\mathcal{Q}$ = set of "strict" precursor payloads
$r$ = circular orbit radius
$\mathcal{S}$ = set of non-linear function samples (Case Study 1), set of vehicles forming a stack (Case Study 2)
$\mathcal{S}'$ = set of vehicle stacks that a vehicle may belong to (Case Study 2)
$t$ = time step
$T$ = maximum number of propellant tanks (Case Study 1)
$x$ = generic design variable (Background and Theory), tank design binary (Case Study 1), commodity flow variable (Case Study 2).
$\beta$ = boil-off rate
$\gamma$ = generated variable
$\Gamma$ = set of generated variables
$\varepsilon$ = binary parameter denoting if an arc is allowed or not (Case Study 2)

$\lambda$ = Lagrange multipliers, or dual variables
$\phi$ = propellant-oxidiser mixture ratio
$\mu$ = gravitational parameter (Case Study 1), ISRU maintenance supply consumption rate (Case Study 2)
$\rho$ = ISRU propellant production rate
$\tau$ = tank mass fraction (Case Study 1), real time-of-flight (Case Study 2)

*Subscripts*
$c$ = commodity type (Case Study 2)
$C, I$ = number of integer commodity types
$C, C$ = number of continuous commodity types
$i$ = generic index (Background and Theory), source node (Case Study 2)
$I$ = number of network nodes (Case Study 2)
$j$ = generic index (Background and Theory), destination node (Case Study 2)
$n$ = manoeuvre index (Case Study 1), logistics vehicle index (Case Study 2)
$o$ = parking orbit index (Case Study 1)
$P$ = number of payloads (Case Study 2)
$R$ = relaxed problem
$s$ = non-linear function sample index
$S$ = set of variables pertaining to a restricted problem
$t$ = tank index (Case Study 1), time step index (Case Study 2)
$v$ = launch vehicle index (Case Study 1)
$V$ = number of logistics vehicles (Case Study 2)
$\oplus$ = Earth

*Superscripts*
cap = capacity
dry = dry mass
$F$ = frequency
$L$ = lower bound
pay = payload mass (Case Study 1), payload capacity (Case Study 2)
prop = propellant
struc = structure
$T$ = transpose
$U$ = upper bound
$*$ = optimal solution
$+$ = inbound flow
$-$ = outbound flow

**Acronyms/Abbreviations**
ACO: Advanced Concepts Office
ConOps: Concept of Operations
KKT: Karush-Kuhn-Tucker conditions
ISRU: *in-situ* resource utilisation
LEO: Low Earth orbit
LP: Linear Programming
MILP: Mixed-Integer Linear Programming
MINLP: Mixed-Integer Non-Linear Programming

**1. Introduction**

Exploration mission architects and program planners frequently face complex sets of categorical decisions when designing the Concepts of Operations (ConOps) of missions. Examples of categorical decisions might be launch vehicle selection, time-sensitive actions such as when to discard a propellant tank, or the choice of a system component design from a discrete set of options.

Approaches for solving complex optimisation problems in space mission design include heuristic optimisation [1], graph theory [2], or linear programming (LP) [3], [4]. Categorical decisions in mission design problems appear as integer variables. Heuristic optimisation is still applicable to such problems [5], [6], but does not guarantee global optimality in a finite amount of computing time. Instead, LP methods can be expanded to operate with integer variables and can still guarantee optimality and give a bound on the objective for suboptimal solutions.

Mixed-integer linear programming (MILP) has been demonstrated as an effective method for optimising large-scale space mission architectures, with prior application to lunar and interplanetary exploration logistics [7], [8], megaconstellation design [9], and on-orbit servicing [10], [11], [12].

Large mixed-integer design problems can become prohibitively computationally expensive to solve with commercially available solvers. In such cases, even if a strong solution can be found, it can be difficult to prove optimality using Strong Duality theory by closing the gap between the solutions of the primal and dual optimisation problems. An example of this is the full model of the second case study of this paper, in which the large lunar logistics model failed to converge to a proven optimal solution. In these situations, approximations can be made in order to restrict the scale of the problem. For example, neural networks can be used to approximate expensive portions of the model [13]. Alternatively, metaheuristics can be used to solve aspects of the problem in combination with LP [14], or used to select appropriate restricted MILP problems in a hierarchical structure [15]. However, these methods can remain computationally expensive due to the need to solve the restricted MILP many times in the search for improved solutions.

Column generation is a method for finding restricted LPs without necessarily requiring additional hierarchical optimization frameworks like metaheuristics. It works by selecting a subset of variables of the problem to form a restricted, with the remaining unselected variables fixed to zero, in a process called the Dantzig-Wolfe decomposition [16]. The unselected variables are then "priced" according to some metric, with the most valuable added back into the restricted variable. If the pricing metric is good, then the objective of the restricted problem will improve with the new set of variables. For continuous variables, the price is the reduced cost of the variable (see Section 2.1). In integer problems and their lack of well-defined dual problem, some heuristic can be employed in order to price the unused variables. This can be as straightforward as solving relaxations of the restricted integer problem, provided that the problem structure doesn't render the relaxation meaningless [17]. Example uses of the Dantzig-Wolfe decomposition applied to the continuous relaxation of MILPs for column generation include multi-depot vehicle scheduling [18] and aircraft scheduling [19]. Wilhelm [20] provides a comprehensive review of further examples of this method.

An alternative approach to solving the pricing problem for integer variables is to construct a dynamic programming subproblem [21], [22]. Both of these methods are somewhat heuristic and do not guarantee that the column generation process produces the restricted MILP that contains the basis of variables that appear in the optimal solution to the full problem.

This paper begins by covering the background and theory of linear programming and duality theory, how this relates to the pricing of variables in the column generation algorithm, and the general space logistics MILP formulation. Next, the paper defines the column generation algorithm used in these studies. Finally, the algorithm is applied to two case studies: first, a human Mars mission involving discrete decision about tank discarding, design, and spacecraft aggregation. The second problem solves a lunar logistics network flow model including discrete decisions about payload launch window assignments. The paper concludes with some discussion of the applicability of column generation to discrete decision making in space mission ConOps design.

## 2. Background and Theory

*2.1 Background on Linear Programming and Duality Theory*

A generic linear optimization problem, or linear program (LP), is shown in Equation 1.

$$\min_{x} f(x) = c^T x$$
$$\text{s.t. } g(x) = Ax - b = 0 \quad (1)$$
$$x \geq 0$$

The set of feasible solutions to an LP forms a polygon bounded by the planes defined by the linear constraints. The vertices of this polygon are called 'basic' solutions. A basic solution $x^*$ can be split into 'basic' and 'non-basic' components, containing zero and non-zero values respectively. The variables in the basic component are called the 'basis'. The number of variables in the basis is equal to the number of constraints in the program.

The Lagrangian for Equation 1 is shown in Equation 2, where $\lambda$ are the Lagrange multipliers.

$$\mathcal{L} = f(x) - \lambda^T g(x) \quad (2)$$

According to Weak Duality theorem, the 'dual' problem to Equation 1 provides a bound on the objective of the original problem, which is called the 'primal' problem. The dual problem consists of a variable for every primal constraint and a constraint for every primal variable. The dual of Equation 1 is shown in Equation 3. Table 1 shows the corresponding dual for all possible components of a primal problem.

$$\max_{\lambda} h(\lambda) = b^T \lambda$$
$$\text{s.t. } A^T \lambda \leq c \quad (3)$$

The Lagrange multipliers $\lambda$ will, from here onwards, be referred to as the dual variables.

According to the Strong Duality theorem, the primal and dual objectives of a linear problem hold the same value when optimality has been reached. Solving the dual of a linear program is extremely useful for understanding how close a particular solution is to optimality. In fact, whether or not a particular solution to the primal problem is optimal can be assessed simply by checking the feasibility of the dual, as for a suboptimal solution, the constraints of Equation 3 are broken. The residual of the dual constraint $c_i - \sum_j A_{i,j} \lambda_j$ is the 'reduced cost', or 'price', of the primal variable $x_i$. The reduced cost of a basic variable will be zero, whilst the reduced cost of a non-basic variable represents the improvement to the objective that can be achieved by moving that variable into the basis.

*2.2 Column Generation*

Column generation is a method for reducing the scale of large linear programs with the aim of improving computation time. It operates by solving the problem using a small subset of variables and leveraging the reduced costs of the unused variables in the resulting sub-optimal solutions. The variables with the greatest reduced costs represent those that would be most useful to add back into the problem.

The process is as follows:
- Select an initial subset of variables that produce a feasible solution.
- Construct and solve the model using this subset of variables. This problem is called the *restricted problem S*.
- Evaluate the full problem *P* using this solution, with the unused variables fixed at 0. Calculate the reduced costs for each unused variable.
- Add the variable with the greatest reduced cost to the restricted problem *S* and repeat.

Table 1. Translation of the components of a linear program between the primal and dual problems.

| | Primal | | Dual |
|---|---|---|---|
| Objective | $\min_x f(x) = c^T x$ | Objective | $\max_\lambda h(\lambda) = b^T \lambda$ |
| Constraint $j$ Sign | $\sum_i A_{i,j} x_i \leq b_j$ | $\lambda_j$ Domain | $\lambda_j \leq 0$ |
| | $\sum_i A_{i,j} x_i \geq b_j$ | | $\lambda_j \geq 0$ |
| | $\sum_i A_{i,j} x_i = b_j$ | | $\lambda_j$ unbounded |
| $x_i$ Domain | $x_i \geq 0$ | Constraint $i$ Sign | $\sum_j A_{i,j} \lambda_j \leq c_i$ |
| | $x_i \leq 0$ | | $\sum_j A_{i,j} \lambda_j \geq c_i$ |
| | $x_i$ unbounded | | $\sum_j A_{i,j} \lambda_j = c_i$ |

The process terminates when there are no remaining unused variables possessing positive reduced costs.

The case studies that follow utilise binary or integer variables to represent categorical design choices in the mission ConOps. The KKT conditions are not meaningful for integer problems because the functions are discontinuous. In this case, the column generation process uses relaxed versions of the problem, $S_R$ being the relaxed restricted problem, and $P_R$ being the relaxation of the full problem $P$. The use of the reduced cost for selecting new variables is therefore a heuristic method when applied to integer variables, because fractional changes to the variables are not feasible in the integer problem. Additionally, variables may have non-zero reduced costs in the relaxed problem, but form part of a degenerate solution once added to the basis, with no change to the relaxed objective. This does not necessarily mean that the added variable will be degenerate in the integer problem. Careful consideration must be made to the structure of the problem as to meaning of fractional changes of otherwise-integer variables and what this means for the reduced costs that will be calculated. The results and discussion sections of this paper will study the heuristic nature of the method.

The overall column generation process for MILPs is shown in Fig. 1, where $x^*_{S,i}$ is the optimal solution to the restricted problem $S_R$ at iteration $i$, and $x_i$ is the same set of values applied to the variables of the full problem $P_R$, with unused variables set to 0.

*2.3 Space Logistics*
The logistics surrounding deep-space exploration missions can be represented using network flow models [4]. Here, vehicles, payloads, and consumables such as propellant are modelled as "commodities" flowing through a graph of nodes (representing locations or specific orbits) connected by arcs (trajectories). In a time-expanded network, the graph is repeated across many time steps, and each repetition is connected by "holdover" arcs that represent the passing of time between the discrete steps.

A network flow model can be represented using a mixed-integer linear program, using the flow of integer or continuous type commodity flows as the decision variables:

$x^-_{i,j,c,t}$: the amount of commodity of type $c$ flowing out of node $i$ to node $j$ at time $t$,
$x^+_{i,j,c,t}$: the amount of commodity of type $c$ flowing into node $j$ from node $i$ at time $t$.

In matrix form, the generic space commodity flow linear program is:

$$\begin{aligned} \min_x \quad & c^T x \\ \text{s.t.} \quad & A_1 x \leq d \\ & A_2 x \leq 0 \\ & A_3 x = 0 \\ & A_4 x = 0 \\ & x \geq 0 \end{aligned}$$

Where $A_1$ is the supply and demand constraint coefficient matrix, $d$ is the supply and demand at each node, $A_2$ is the commodity capacity constraint coefficient matrix, $A_3$ is the commodity dynamics constraint coefficient matrix, and $A_4$ is a constraint coefficient matrix that limits commodities to only the allowed set of arcs. Case Study 2 of this papers concerns a lunar exploration logistics scenario, and provides the specific formulations of each of these constraints for the commodity types, nodes, and arcs that are relevant to that scenario.

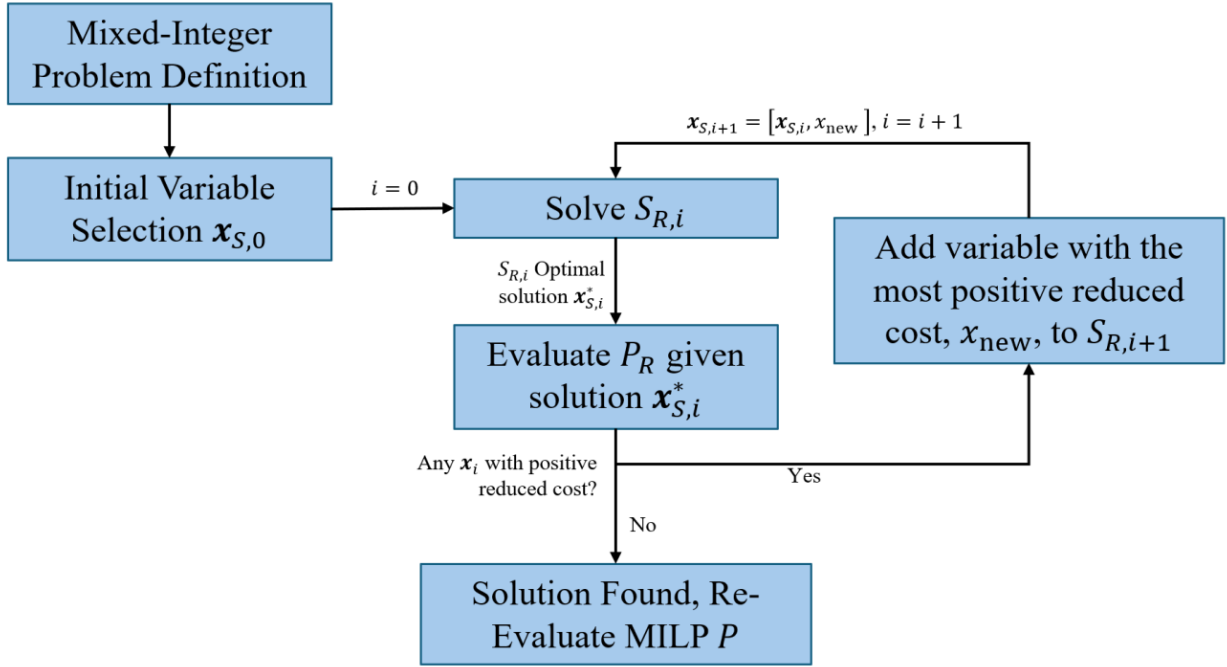

Fig. 1: Flow chart of the column generation process.

### 3. Method
*3.1 Column Generation Problem Ontology*

It is not necessary to remove *all* variables from the problem before initiating the column generation process. Removing specific sets of variables can improve the solve time of the problem more than others, and generating across a subset of variables results in a smaller pricing problem. Therefore, some attention must be given to the selection of the set of variables to be generated. Here, the terminology and notation surrounding this choice is explained.

First, the set of variables to be removed from the problem and then generated is referred to as $\Gamma$. The column generation algorithm could generate variables individually, although this may be a slow process with little or no change to the optimal solution per column generation iteration. Alternatively, variables can be generated collectively. At each iteration of the column generation algorithm, the group with the best price is the one that is generated. The price of variable group $\Gamma_i$ is $\pi_i$:

$$\pi_i = \sum_{y \in \Gamma_i} \pi_y = \sum_{y \in \Gamma_i} \left( c_y - \sum_k A_{y,k} \lambda_k \right) \quad (4)$$

Where $\pi_y$ is the price of variable $y$. Note that the process can be terminated before all groupings $\Gamma_i \in \Gamma$ have been generated if there is no expected improvement to the objective, or if other termination conditions are reached.

*3.2 Strategy for Defining Variable Groups*

Various strategies could be taken for defining variables groupings – this section describes the strategy taken in this paper and is specific to the problem structures involved with the case studies. The case studies presented this paper involve mixed-integer linear programs containing several categories of variables, which may be continuous or integer, and each have their own sets of indices, e.g. the variable set $x_{t,n,s,v}$ from Case Study 1 has indices $t, n, s,$ and $v$. If $x$ is chosen as the set of variables to be initially removed from the problem ($\Gamma$), and then re-added through column generation, a subset of its index set $(t, n, s, v)$ is defined as a *grouping* index set. So, a generating variable group $\Gamma_i$ is comprised of variables that share the same grouping index, and represent all possible combinations of the remaining indices.

In general terms: let $l$ be the number of dimensions of a variable set and $\mathcal{I}_d$ be the set of values to be included from index $d$. $\mathcal{I}_1 \times \mathcal{I}_2 \times ... \times \mathcal{I}_l$ is then the topological space spanned by the variable set. Continuing the above example, if $x$ is a four-index variable with $\mathcal{I}_t = \{0,1,2,...,11\}, \mathcal{I}_n = \{0,1,2,3,4\}, \mathcal{I}_s = \{0,1,2,3,...,29\},$ and $\mathcal{I}_v = \{0,1,2,...,69\}$, then $x \in \mathbb{Z}^{12 \times 5 \times 30 \times 70}$.

Let $\mathcal{I}^G$ be the grouping index set, whose cardinality is $m$. Let the set $G$ be the set of all possible combinations of the grouping indices:

$$G = \{(i_1, i_2, \ldots, i_m) \ \forall \ i_1 \in \mathcal{I}_1^G, i_2 \in \mathcal{I}_2^G, \ldots, i_m \in \mathcal{I}_m^G\}$$

Continuing the example three-dimensional $x \in \mathbb{Z}^{12 \times 5 \times 30 \times 70}$, if $\mathcal{I}^G = \{\mathcal{I}_s, \mathcal{I}_v\}$, then

$$G = \{(0,0), (0,1), \ldots, (0,69), (1,0), (1,1), \ldots, (1,69), \ldots, (29,0), (29,1), \ldots, (29,69)\}$$

The set of remaining indices, all combinations of which are included in each generating index set, is then $\mathcal{I}^N = \mathcal{I} \setminus \mathcal{I}^G$:

$$N = \{\{(j_1, j_2, \ldots, j_{l-m}) \ \forall \ j_1 \in \mathcal{I}_1^N, j_2 \in \mathcal{I}_2^N, \ldots, j_{l-m} \in \mathcal{I}_{l-m}^N\}$$

Continuing our example, $\mathcal{I}^N = \{\mathcal{I}_t, \mathcal{I}_n, \mathcal{I}_s, \mathcal{I}_v\} \setminus \{\mathcal{I}_s, \mathcal{I}_v\} = \{\mathcal{I}_t, \mathcal{I}_n\}$, and:

$$N = \{(0,0), (0,1), \ldots, (0,4), (1,0), (1,1), \ldots, (1,4), \ldots, (11,0), (11,1), \ldots, (11,4)\}$$

For example, the generating group $\Gamma_i$ for the generating index $i = (0,0)$ is $\Gamma_{0,0} = \{x_{0,0,0,0}, \ldots, x_{11,4,0,0}\}$.

The price of an individual group of generated variables is:

$$\pi_i = \sum_{j \in N} \pi_{i,j} = \sum_{j \in N} \left( c_{i,j} - \sum_k A_{i,j,k} \lambda_k \right) \quad (5)$$

Where $\pi_{i,j}$ is the price of variable $x_{i,j}$. If there are multiple classes of variables in Γ with different index dimensions, each variable must have its own $G$ and $N$ defined, and the prices are calculated separately. The user can decide whether to take the most valuable group overall, or the most valuable group for each variable.

*3.3 Selection of Initial Variable Set*

The selection of the initial variable set must provide a feasible starting point for the column generation process. This was done by one of two methods, both involving the construction of a separate binary optimisation problem. The first method, shown in Equation 6, minimises the number of generated variables required to satisfy the constraints that contain *only* generated variables.

$$\min_\gamma \sum_i \gamma_i$$
$$\text{s.t.} \ g(\gamma) = b \ \forall \ \gamma \in \Gamma \quad (6)$$
$$\gamma \in \{0,1\}$$

Where $\gamma$ is the set of non-zero variables in $\Gamma$. In cases where the constraints containing both generated and non-generated variables are not overly complex, this method is effective at producing initial starting points for the column generation process.

Otherwise, a problem can be solved with the same objective to minimise the number of non-zero generated variables, but with all constraints from the full problem *P* included. This is shown in Equation 7.

$$\min_\gamma \sum_i \gamma_i$$
$$\text{s.t.} \ g(\gamma, x) = b \ \forall \ \gamma \in \Gamma, x \notin \Gamma \quad (7)$$
$$\gamma \in \{0,1\}, x \geq 0$$

Equation 7 is slower to solve than Equation 6, but still typically faster than the full problem *P* due to the simpler, integer objective.

**4. Case Study Setup**

*4.1 Case Study 1: Crewed Mars Mission*

The first case study concerns a crewed exploration mission to Mars [23]. The ConOps begin with the assembly of the spacecraft in low Earth orbit (LEO). The spacecraft then injects itself onto a trans-Martian trajectory, and captures into Mars orbit some months later. The crew carries out their mission to the surface and then returns to the same spacecraft, which remained in orbit. The spacecraft then injects into a trans-Earth orbit, with two deep-space corrections en-route. The spacecraft then captures into LEO, with the crew subsequently returning to Earth via a separate re-entry vehicle. An example ConOps is shown in Fig. 2, though the referenced study assumed spacecraft assembly in MEO and lunar NRHO rather than in LEO. This difference does not affect the methodology discussed in this paper.

In this simplified model, the spacecraft is assumed to consist of a structure, a habitat, and up to 12 propellant tanks. The tanks can be dropped once they are expended in order to save mass. The spacecraft is assumed to use nuclear thermal propulsion, with an engine with an $I_{sp}$ of 900 s. The payload mass carried by the spacecraft varies throughout the mission, representing the use of crew consumables throughout the interplanetary transit, and during the crew's time on Mars.

The aims of this case study were to find the following:
1. optimal tank drop sequence,
2. tank sizing,
3. propellant allocation,
4. launch vehicle usage,
5. and initial parking orbit for assembly.

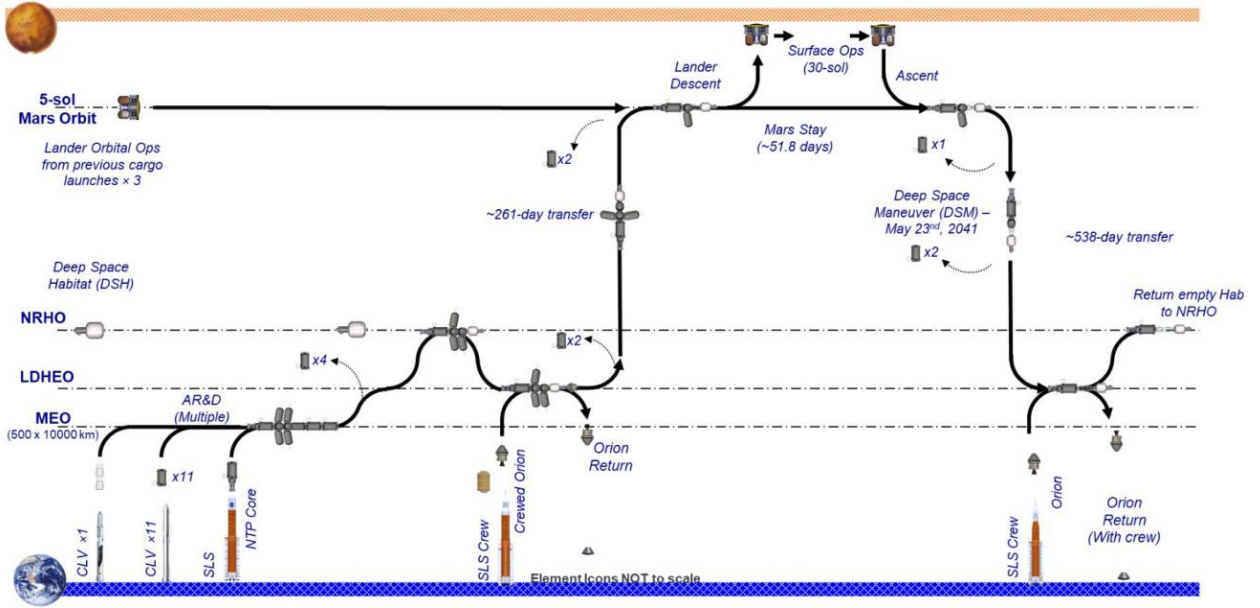

Fig. 2: Example of a crewed Mars mission Concept of Operations featuring tank drop staging. Adapted from [25].

Of these 5 points, items 1, 4, and 5 are categorical decisions and will be represented by binary variables. 3 is strictly a continuous variable (confirming this to be a mixed-integer problem), whilst item 2 can be modelled in a number of ways. A simplistic method would be to assume that the tank dry mass is linearly proportional to its propellant capacity, using some rule-of-thumb percentage. A higher fidelity approach would be to develop a non-linear relationship between tank dry mass and capacity by studying real-world data, such as the data shown in Fig. [24]. The relationship between the tank dry mass-to-propellant capacity ratio (from here onwards referred to as tank fraction $\tau$), and propellant capacity $m^{\text{cap}}$ is indicated by the parametric function shown in Equation 8.

$$\tau = 0.7699 \, (m^{\text{cap}})^{-0.187} \qquad (8)$$

Implementing this relationship directly in the MILP would create a mixed-integer non-linear program (MINLP), a problem type which is notoriously difficult to solve and is not the subject of this paper (see [26] for how such problems can be tackled).

Given the presence of binary variables relating to the propellant tanks already in this model, an alternative approach would be to take as set discrete samples $\mathcal{S}$ from Equation 8 over the expected range of the tank capacities, and add an index to the tank drop sequence binary variable that determines which sample to choose for the tank sizing. In this case study, 30 linearly spaced samples of Equation 8 were taken over the range of $m^{\text{cap}} \in [5000, 40000]$ kg.

Launch vehicle and parking orbit selection were also a categorical choice to be made using binary variables. Both of these decisions combine to determine the mass limit of the individual spacecraft components, because a specific launch vehicle has a maximum payload capacity that is specific to its target parking orbit. To formulate these mass constraints, the NASA Launch Vehicle Performance website was consulted [27]. The payload capability of a variety of U.S. launch vehicles versus 28.5° LEO altitude is shown in Fig. 4. The total number of launches that a payload can be assigned to is equal to the product of the number of available launch vehicle *types* and the number of items to be launched. This encompasses all payload-to-launch vehicle packing options. For example, if there are 5 launch vehicles types and 10 payloads, then there can be up to 10 launches of each vehicle type, for a total number of 50 launch options.

These altitudes were converted to $\Delta V$ requirements using Equation 9 by assuming that the spacecraft would be injecting itself into a trans-Mars orbit with a $C_3$ of 16 km$^2$/s$^2$ [28]. $r$ is the radius of the parking orbit in question, assuming the parking orbit to be circular.

$$\Delta V = \sqrt{C_3 + \frac{2\mu_\oplus}{r}} - \sqrt{\frac{\mu_\oplus}{r}} \qquad (9)$$

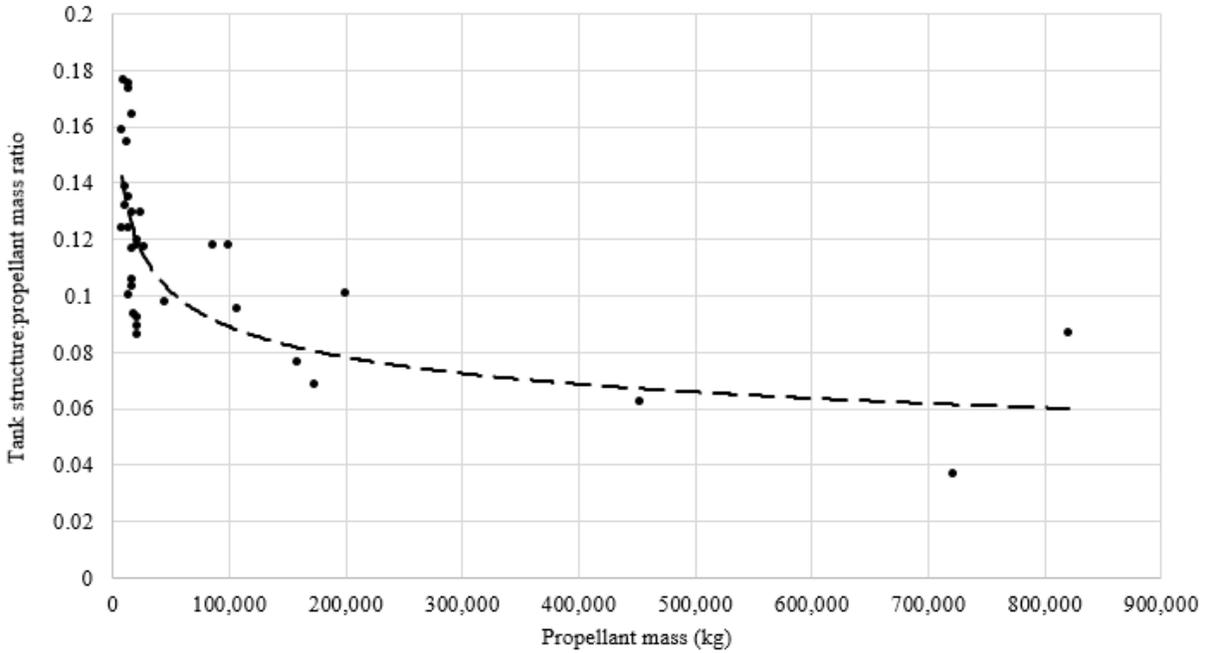

Fig. 3: Parametric data showing propellant capacity versus dry mass-to-propellant capacity ratio for a range of real-world propellant tanks.

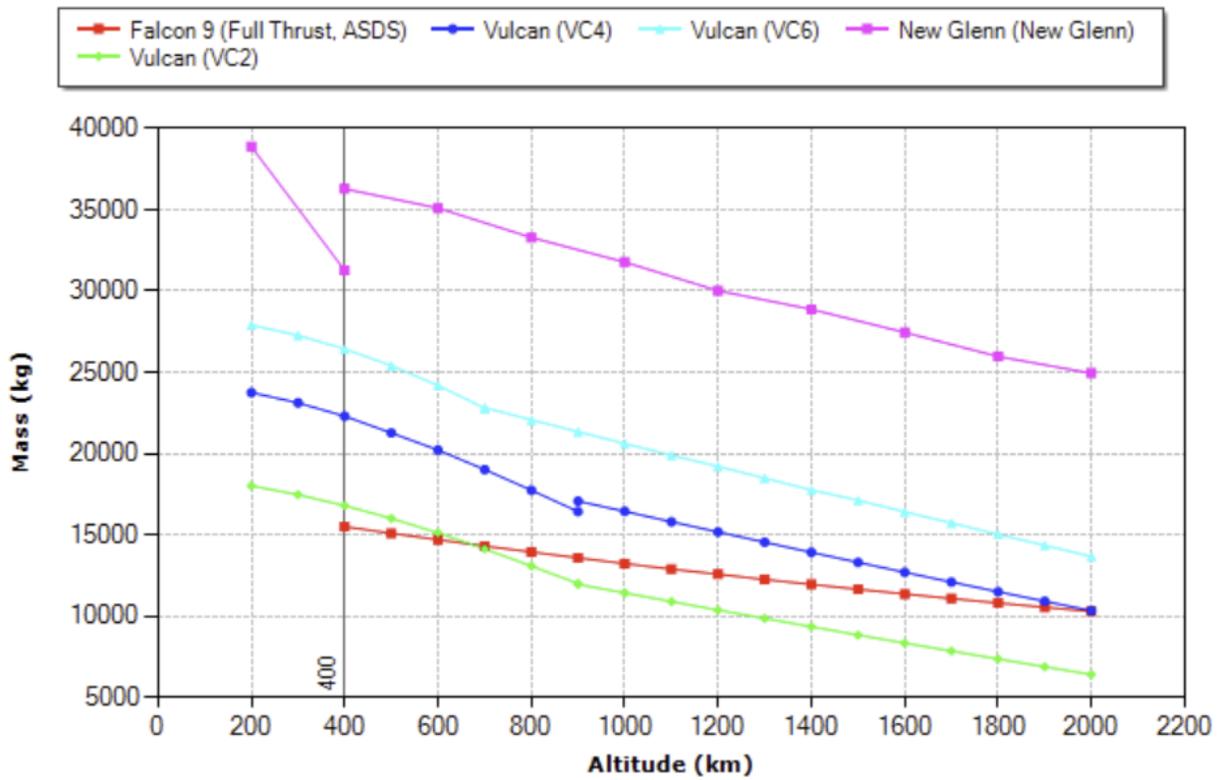

Fig. 4: Payload capability versus parking orbit altitudes, with inclination 28.5°, for a variety of launch vehicles.

IAC-24-D2,8,x85968 Page 8 of 21

The problem variables are:

- $x_{t,n,s,v}$: binary variable determining if tank $t$ of design $s$ is present during manoeuvre $n$, and was launched by launch vehicle $v$.
- $m_{t,n,v}^{\text{prop}}$: continuous variable determining the amount of propellant drawn from tank $t,v$ during manoeuvre $n$.
- $y_v^{\text{struc}}$: binary variable determining the launch vehicle $v$ that the spacecraft structure is launched on.
- $y_v^{\text{pay}}$: binary variable determining the launch vehicle $v$ that the initial payloads (habitat, crew and supplies) are launched on.
- $\ell_{v,o}$: binary variable determining whether spacecraft $v$ is launched into parking orbit $o$. Note that repeat launches of the same vehicle are allowed, so $|\mathcal{V}|$ is the number of launcher types multiplied by the number of payloads. This allows any combination of payload-to-vehicle assignment.
- $\sigma_o$: binary variable determining if parking orbit $o$ is selected.

Combining all of these variables and constraints into a single MILP, the problem formulation is as follows.

$$\min_{m^{\text{prop}},x,y,\sigma,\ell} \sum_v \sum_t \left( \sum_n m_{t,n,v}^{\text{prop}} + \sum_s x_{t,n,s,v} m_s^{\text{cap}} \tau_s \right) \quad (10)$$

s.t. 1A: $\sum_v \sum_t m_{t,n=0,v}^{\text{prop}}$
$$\geq -\mathcal{M}(1-\sigma_o)$$
$$+ (Z_o - 1)\left( m^{\text{struc}} + m_0^{\text{pay}} \right.$$
$$+ \sum_v \sum_t \left( \sum_s x_{t,n=0,s,v} m_s^{\text{cap}} \tau_s \right.$$
$$\left.\left. + \sum_{i=1}^N m_{t,i,v}^{\text{prop}} \right) \right) \; \forall \, o \in \mathcal{O} \quad (11)$$

1B: $\sum_v \sum_t m_{t,n,v}^{\text{prop}} \geq (Z_n - 1)\left( m^{\text{struc}} + \sum_{i=0}^n m_i^{\text{pay}} \right.$
$$+ \sum_v \sum_t \left( \sum_s x_{t,n,s,v} m_s^{\text{cap}} \tau_s \right.$$
$$\left.\left. + \sum_{i=n+1}^N m_{t,i,v}^{\text{prop}} \right) \right) \; \forall \, n \in \{1 \ldots N\}$$
$$(12)$$

2A: $\sum_n m_{t,n,v}^{\text{prop}} \leq \sum_s x_{t,n=0,s,v} m_s^{\text{cap}} \; \forall \, t \in T, v \in \mathcal{V}$ (13)

2B: $m_{t,n,v}^{\text{prop}} \leq \sum_s x_{t,n,s,v} m_s^{\text{cap}} \; \forall \, t \in T, n \in N, v \in \mathcal{V}$ (14)

3A: $x_{t,i,s,v} \leq x_{t,n,s,v} \; \forall \, t \in T, n \in N, s \in \mathcal{S}, v \in \mathcal{V}, i \in \{n \ldots N\}$
$$(15)$$

3B: $\sum_v \sum_s x_{t=T,n=N,s,v} = 1$ (16)

4A: $\sum_s \sum_v x_{t,n=0,s,v} \leq 1 \; \forall \, t \in \{0,\ldots,T\}$ (17)

4B: $\sum_v y_v^{\text{struc}} = 1$ (18)

4C: $\sum_v y_v^{\text{pay}} = 1$ (19)

5: $\sum_t \left( \sum_s x_{t,n=0,s,v} m_s^{\text{cap}} \tau_s + \sum_n m_{t,n,v}^{\text{prop}} \right)$
$$+ y_v^{\text{struc}} m^{\text{struc}} + y_v^{\text{pay}} m_0^{\text{pay}}$$
$$\leq \sum_o m_{v,o}^{LV} \ell_{v,o} \; \forall \, v \in \mathcal{V}$$
$$(20)$$

6A: $\ell_{v,o} \leq \sigma_o \; \forall \, v \in \mathcal{V}, o \in \mathcal{O}$ (21)

6B: $\sum_o \sigma_o = 1$ (22)

The objective of the MILP is to minimize the total tank and propellant mass launched into the parking orbit. The constraint explanations are as follows:

1. A: Propellant consumption during the first manoeuvre obeys the rocket equation, with mass fraction $Z_o$ calculated according to the trans-Mars injection $\Delta V$ associated with parking orbit $o$. The parameter $\mathcal{M}$ is a sufficiently large coefficient such that this constraint cannot be activated for orbits that are not selected. This formulation is intended to avoid any bi-linear constraints.
   B: Propellant consumption for subsequent manoeuvres obeys the rocket equation with standardized $\Delta V$, regardless of initial orbit selection. These $\Delta V$'s are listed in Table 2. The changes in payload mass after each part of the mission are also listed in Table 2.
2. A: Total propellant across all manoeuvres assigned to a tank must be within that tank's capacity.
   B: Propellant consumed within a specific manoeuvre must be within the tanks capacity (prevents propellant consumption from tanks that are no longer present).
3. A: Once a tank has been dropped, it cannot be re-attached.
   B: The final tank remains attached for the entire mission.
4. A: Select only one launch vehicle / design sample pair for each tank.
   B: Only select one launch vehicle for the spacecraft structure.
   C: Only select one launch vehicle for the payload.
5. The total mass assigned to a single launch vehicle must be within the launch vehicle's capacity.
6. A: Only launch vehicles into the selected parking orbit.
   B: Choose exactly one parking orbit.

Table 2: $\Delta V$'s used throughout the Mars mission ConOps, based on the 2039 opposition mission from [29].

| # | Manoeuvre | $\Delta V$ (m/s) | $\Delta m^{\text{pay}}$ (kg) |
|---|---|---|---|
| 0 | Trans-Mars Injection | $Z_o$ | 25000 |
| 1 | Mars capture | 1429 | -500 |
| 2 | Trans-Earth injection | 1908 | -20000 |
| 4 | DSM 2 | 0 | -500 |
| 5 | Earth capture | 4500 | -500 |

The crewed Mars mission ConOps are optimised using the launch vehicle variable $\ell_{v,o}$ as the set of variables to be generated, using the vehicle index $v$ as the grouping index. The initial set of variable indices for the launch vehicle were generated using Equation 7, and included [58, :], [60, :], [61, :], [62, :], [64, :], [66, :], and [69, :], with ':' indicating that all indices in the other dimensions are available to the problem. These all correspond to New Glenn vehicles.

*4.2 Case Study 2: Lunar Logistics Model*
The second case study analyses the logistics of a long-term lunar exploration program, based on the future phases of the Artemis program as laid out by the Global Exploration Roadmap [30], [31], [32]. The model developed in this case study finds the optimal flow of goods through a logistics network consisting of the Earth, Low Lunar Orbit, and the Moon's surface, such that supply and demand constraints defined by exploration program requirements are fulfilled. Further constraints are imposed based on logistics vehicle availability. The logistics network is shown in Figure , adapted from [33], and the costs associated with each arc are listed in Table 4. The arcs between different nodes represent the trajectories between them, and the arcs between nodes and their future counterparts represent stays within that location. Each pair of time steps forms one lunar period ($\approx$ 29.53 days). This is referred to as a macro-period, and each step within the macro-period as a micro-period. The length of the holdover arc between the micro-periods *within* a macro-period is consistent with a three-day mission on the lunar surface, in line with the later Apollo missions. The holdover arc connecting the second micro-period with the first of the *next* macro-period has a "time of flight" such that the total length of the macro-period is made up to a lunar period.

The flow of goods is restricted based on the current time index. Outbound (Earth towards Moon) flow is only allowed on even time indices, whilst return (Moon towards Earth) flow is only allowed on odd time indices. This prevents any opposing commodity flows from "cancelling each other out" and resulting in an unrealistic net loss of goods.

The types of commodities considered are listed in Table 3.

Table 3: Commodity types considered in the lunar logistics model and their unit masses.

| $c$ | Name | Variable Type | $m_{n,c}$ (kg) |
|---|---|---|---|
| 0 | Vehicle | Integer | $m_n^{\text{dry}}$ |
| 1 | Crew | Integer | 100 |
| 2 | ISRU infrastructure | Continuous | 1 |
| 3 | Maintenance Supplies | Continuous | 1 |
| 4 | Crew Consumables | Continuous | 1 |
| 5 | Inert Payload | Continuous | 1 |
| 6 | Propellant Oxidiser | Continuous | 1 |
| 7 | Propellant Fuel | Continuous | 1 |

$I_{\text{sp}}$ was used as a proxy for propellant type. Vehicle with $I_{\text{sp}} < 370$ s were assumed to have storable propellant, $370 \leq I_{\text{sp}} < 420$ s were assumed to have liquid methane/liquid oxygen propellant, and those with $I_{\text{sp}} > 420$ s were assumed to have liquid hydrogen/liquid oxygen propellant. Each propellant types had different boil-off rates associated with them.

Additionally, the vehicles that take part in the campaign are allowed to form vehicle "stacks", which is a system to model the rendezvous and separation of vehicles at nodes (but not between them). The method is based on [34]. Vehicle stacks are formulated as a list of their constituent vehicle, taking the combined dry mass and propellant capacities, $I_{\text{sp}}$ of the lead vehicle, and payload capacity of the lead vehicle minus the dry masses of the carried vehicles.

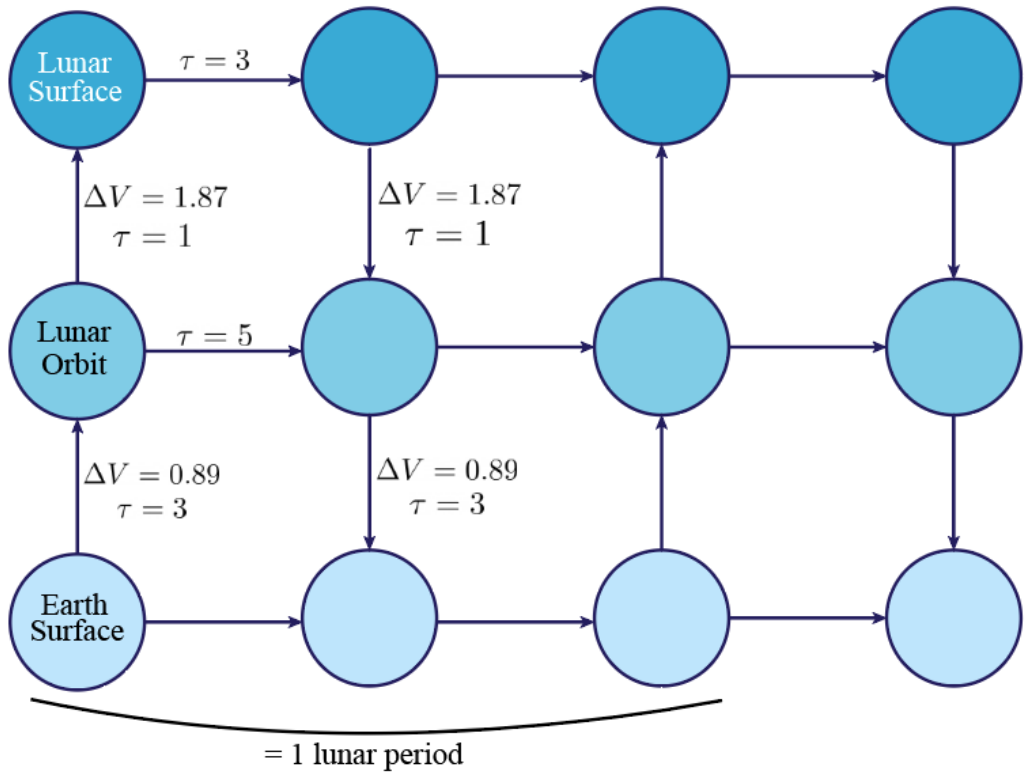

Figure 5: Time-expanded network of nodes (0 = Earth. 1 = Low Lunar Orbit, 2 = Lunar Surface, connected by arcs. "Holdover" arcs representing the passing of discrete time steps.

Table 4: $\Delta V$ and time of flight of each arc in the logistics network. Real time of flight measured in ° sun angle, equivalent of 1/360th of a Sun-Moon-Earth synodic period in the bicircular restricted 4 body problem.

| Arc $[i,j]$ | Launch cost | $\Delta V$ (km/s) | Real TOF (° sun angle) | Discrete TOF |
|---|---|---|---|---|
| [0,0] | - | 0 | - | 1 |
| [0,1] | 1 | 0.89 | 36.57 | 0 |
| [1,1] | - | 0.15 per year [ref] | 60.95 on even time step, 360-60.95 otherwise | 1 |
| [1,2] | - | 1.87 | 12.19 | 0 |
| [2,2] | - | - | 36.57 on even time step, 360-36.57 otherwise | 1 |
| [2,1] | - | 1.87 | 12.19 | 0 |
| [1,0] | - | 0.89 | 36.57 | 0 |

The lunar logistics model was formulated as the following MILP, using the following variables:

- $x^+_{n,i,j,c,t}$: flow of goods of commodity type $c$ from node $i$ towards node $j$, carried by vehicle $n$ at time $t$. Can be integer or continuous depending on $c$.
- $x^-_{n,i,j,c,t}$: flow of goods of commodity type $c$ from node $i$ arriving at node $j$, carried by vehicle $n$ at time $t$. Can be integer or continuous depending on $c$.
- $B_{p,t}$: binary variable determining the demand time index $t$ of payload $p$.

$$\min_{x,B} f(x) = \sum_t \sum_n \sum_c m_{n,c}\, x^-_{n,i=0,j=1,c,t} \tag{23}$$

$$\text{s.t. 1A:} \sum_{n' \in \mathcal{S}'_n} \sum_j \left( x^-_{n,i,j,c=0,t} - x^+_{n,j,i,c=0,t} \right) \le \sum_n d_{n,i,c=0,t} \quad \forall\, n, i, t \tag{24}$$

$$\text{1B:} \sum_n \sum_j \left( x^-_{n,i,j,c,t} - x^+_{n,j,i,c,t} \right) \le \sum_n \left( d^f_{n,i,c,t} + \sum_p B_{p,t} d^m_{n,i,c,p} \right) \quad \forall\, i, t, c > 0 \tag{25}$$

$$\text{2A:} \sum_{c \in \{1,2,3,4,5\}} \left( m_{n,c} x^-_{n,i,j,c,t} \right) - m^{\text{pay}}_n x^-_{n,i,j,c=0,t} \le 0 \quad \forall\, n, i: i \ne 2, j: j \ne 2, t \tag{26}$$

$$\text{2B:} \sum_{c \in \{1,3,4,5\}} \left( m_{n,c} x^-_{n,i=2,j=2,c,t} \right) - m^{\text{pay}}_n x^-_{n,i=2,j=2,c=0,t} \le 0 \quad \forall\, n, t \tag{27}$$

$$\text{2C:}\ x^-_{n,i,j,c=6,t} - \phi_n m^{\text{prop}}_n x^-_{n,i,j,c=0,t} \le 0 \quad \forall\, n, i, j, t \tag{28}$$

$$\text{2D:}\ x^-_{n,i,j,c=7,t} - (1-\phi_n) m^{\text{prop}}_n x^-_{n,i,j,c=0,t} \le 0 \quad \forall\, n, i, j, t \tag{29}$$

$$\text{3A:}\ x^+_{n,i,j,c=4,t} - x^-_{n,i,j,c=4,t} + k\tau_{n,i,j,t} x_{n,i,j,c=1,t} = 0 \quad \forall\, n, i, j, t \tag{30}$$

$$\text{3B:} \sum_n \left( x^+_{n,i=2,j=2,c=3,t} - x^-_{n,i=2,j=2,c=3,t} + \mu\tau_{n,i=2,j=2,t} x_{n,i=2,j=2,c=2,t} \right) = 0 \quad \forall\, t \tag{31}$$

$$\text{3C:}\ x^+_{n,i,j,c=3,t} - x^-_{n,i,j,c=3,t} = 0 \quad \forall\, n, t, i: i \ne 2, j: j \ne 2 \tag{32}$$

$$\text{3D:}\ x^+_{n,i,j,c=6,t} - (1-\beta^{\text{ox}}_n)^{\tau_{n,i,j,t}} x^-_{n,i,j,c=6,t} + \phi Z_{n,i,j,t} \left( \sum_c m_{n,c} x^-_{n,i,j,c,t} \right) = 0 \quad \forall\, n, t, (i,j): \{i \ne j\} \cup \{i = j \ne 2\} \tag{33}$$

$$\text{3E:}\ x^+_{n,i=2,j=2,c=6,t} - (1-\beta^{\text{ox}}_n)^{\tau_{n,i=2,j=2,t}} x^-_{n,i=2,j=2,c=6,t} - \rho\tau_{n,i=2,j=2,t} x^-_{n,i=2,j=2,c=2,t} = 0 \quad \forall\, n, t \tag{34}$$

Table 5: Commodity types considered in the lunar logistics model and their unit masses [include refs].

| $n$ | Name | $m^{\text{pay}}_n$ | $m^{\text{prop}}_n$ | $m^{\text{dry}}_n$ | $I_{\text{sp}}$ | $t^f$ | $t^L$ |
|---|---|---|---|---|---|---|---|
| 0 | Astrobotic Griffin | 625 | 3320 | 1950 | 340 | 12 | 0 |
| 1 | B.O. Blue Moon | 4500 | 6350 | 2150 | 420 | 6 | 0 |
| 2 | Draper/ispace S2 | 500 | 3380 | 2120 | 340 | 12 | 0 |
| 3 | I.M. Nova-C | 100 | 1010 | 790 | 370 | 6 | 0 |
| 4 | ESA EL3 | 1800 | 5580 | 2520 | 340 | 36 | 0 |
| 5 | ISECG lander | 9000 | 23660 | 9340 | 340 | 12 | 0 |
| 6 | ISECG ascender | 500 | 10000 | 1000 | 340 | 12 | 0 |
| 7 | Orion | 11800 | 22000 | 16520 | 316 | 1 | 0 |
| 8 | MK2 ISECG lander | 11390 | 23660 | 9340 | 370 | 12 | 0 |
| 9 | MK2 ISECG ascender | 500 | 10000 | 1000 | 370 | 12 | 0 |

Table 6: Payloads to be launched in the Artemis campaign, including quantities, launch windows, and precursor or co-payload requirements.

| $p$ | Name | $c$ | Quantity | $i_p$ | $j_p$ | $t_p^L$ | $t_p^U$ | $\mathcal{P}_p$ | $\mathcal{Q}_p$ | $\mathcal{C}_p$ |
|---|---|---|---|---|---|---|---|---|---|---|
| 0 | Power Plant Element | 5 | 1500 | 0 | 2 | 0 | 48 | | | |
| 1 | Artemis 7 Crew | 1 | 4 | 0 | 1 | 0 | 48 | 0 | | |
| 2 | Artemis 7 Crew Landing | 1 | 4 | 1 | 2 | 0 | 48 | | | 1 |
| 3 | Artemis 7 Crew Ascent | 1 | 4 | 2 | 1 | 0 | 48 | | 2 | |
| 4 | Artemis 7 Crew Return | 1 | 4 | 1 | 0 | 0 | 48 | | | 3 |
| 5 | Sample Return | 5 | 200 | 2 | 0 | 0 | 48 | | | |
| 6 | Habitat | 5 | 4500 | 0 | 2 | 0 | 54 | | | |
| 7 | Artemis 8 Crew | 1 | 4 | 0 | 1 | 0 | 54 | 6 | 4 | |
| 8 | Artemis 8 Crew Landing | 1 | 4 | 1 | 2 | 0 | 54 | | | 7 |
| 9 | Artemis 8 Crew Ascent | 1 | 4 | 2 | 1 | 0 | 54 | | 8 | |
| 10 | Artemis 8 Crew Return | 1 | 4 | 1 | 0 | 0 | 54 | | | 9 |
| 11 | Sample Return | 5 | 200 | 2 | 0 | 0 | 54 | | | |
| 12 | Artemis 9 Crew | 1 | 4 | 0 | 1 | 12 | 60 | | 10 | |
| 13 | Artemis 9 Crew Landing | 1 | 4 | 1 | 2 | 12 | 60 | | | 12 |
| 14 | Artemis 9 Crew Ascent | 1 | 4 | 2 | 1 | 12 | 60 | | 13 | |
| 15 | Artemis 9 Crew Return | 1 | 4 | 1 | 0 | 12 | 60 | | | 14 |
| 16 | Sample Return | 5 | 200 | 2 | 0 | 12 | 60 | | | |
| 17 | Pressurised Rover | 5 | 4500 | 0 | 2 | 0 | 66 | | | |
| 18 | Pressurised Rover | 5 | 4500 | 0 | 2 | 0 | 66 | | | |
| 19 | Artemis 10 Crew | 1 | 4 | 0 | 1 | 24 | 66 | 17, 18 | 15 | |
| 20 | Artemis 10 Crew Landing | 1 | 4 | 1 | 2 | 24 | 66 | | | 19 |
| 21 | Artemis 10 Crew Ascent | 1 | 4 | 2 | 1 | 24 | 66 | | 20 | |
| 22 | Artemis 10 Crew Return | 1 | 4 | 1 | 0 | 24 | 66 | | | 21 |
| 23 | Sample Return | 5 | 200 | 2 | 0 | 24 | 66 | | | |
| 24 | Artemis 11 Crew | 1 | 4 | 0 | 1 | 36 | 72 | | 22 | |
| 25 | Artemis 11 Crew Landing | 1 | 4 | 1 | 2 | 36 | 72 | | | 24 |
| 26 | Artemis 11 Crew Ascent | 1 | 4 | 2 | 1 | 36 | 72 | | 25 | |
| 27 | Artemis 11 Crew Return | 1 | 4 | 1 | 0 | 36 | 72 | | | 26 |
| 28 | Sample Return | 5 | 200 | 2 | 0 | 36 | 72 | | | |
| 29 | Artemis 12 Crew | 1 | 4 | 0 | 1 | 48 | 84 | | 27 | |
| 30 | Artemis 12 Crew Landing | 1 | 4 | 1 | 2 | 48 | 84 | | | 29 |
| 31 | Artemis 12 Crew Ascent | 1 | 4 | 2 | 1 | 48 | 84 | | 30 | |
| 32 | Artemis 12 Crew Return | 1 | 4 | 1 | 0 | 48 | 84 | | | 31 |
| 33 | Sample Return | 5 | 200 | 2 | 0 | 48 | 84 | | | |
| 34 | Fission Power Plant | 5 | 4500 | 0 | 2 | 48 | 84 | | | |
| 35 | Habitat | 5 | 4500 | 0 | 2 | 48 | 84 | | | |
| 36 | Artemis 13 Crew | 1 | 4 | 0 | 1 | 60 | 96 | 34, 35 | 32 | |
| 37 | Artemis 13 Crew Landing | 1 | 4 | 1 | 2 | 60 | 96 | | | 36 |
| 38 | Artemis 13 Crew Ascent | 1 | 4 | 2 | 1 | 60 | 96 | | 37 | |
| 39 | Artemis 13 Crew Return | 1 | 4 | 1 | 0 | 60 | 96 | | | 38 |
| 40 | Sample Return | 5 | 200 | 2 | 0 | 60 | 96 | | | |
| 41 | Artemis 14 Crew | 1 | 4 | 0 | 1 | 72 | 96 | | 39 | |
| 42 | Artemis 14 Crew Landing | 1 | 4 | 1 | 2 | 72 | 96 | | | 41 |
| 43 | Artemis 14 Crew Ascent | 1 | 4 | 2 | 1 | 72 | 96 | | 42 | |
| 44 | Artemis 14 Crew Return | 1 | 4 | 1 | 0 | 72 | 96 | | | 43 |
| 45 | Sample Return | 5 | 200 | 2 | 0 | 72 | 96 | | | |

$$3F: x^+_{n,i,j,c=7,t} - \left(1 - \beta_n^f\right)^{\tau_{n,i,j,t}} x^-_{n,i,j,c=7,t}$$
$$+ (1 - \phi_n) Z_{n,i,j,t} \left(\sum_c m_{n,c} x^-_{n,i,j,c,t}\right)$$
$$= 0 \; \forall \; n, i, j, t \quad (35)$$

$$3G: x^+_{n,i,j,c,t} - x^-_{n,i,j,c,t} = 0 \; \forall \; n, i, j, t, c \in \{0,1,5\} \quad (36)$$

$$4: x^-_{n,i,j,c,t} = 0 \; \forall \; c, (n,i,j,t): \varepsilon_{n,i,j,t} = 0 \quad (37)$$

$$5A: \sum_t B_{p,t} = 1 \; \forall \; p \quad (38)$$

$$5B: B_{p,t} = B_{p',t} \; \forall \; t, p, p' \in \mathcal{C}_p \quad (39)$$

$$5C: \sum_{t'}^{t} B_{p,t'} \leq \sum_{t'}^{t} B_{p',t'} \leq t, p, p' \in \mathcal{P}_p \quad (40)$$

$$5D: \sum_{t'}^{t} B_{p,t'} \leq \sum_{t'}^{t-2} B_{p',t'} \leq t, p, p' \in \mathcal{Q}_p \quad (41)$$

$$5E: B_{p,t} = 0 \; \forall \; p, t: (t < t_p^L) \text{ or } (t > t_p^U) \text{ or } (i_p > j_p \text{ and } t \text{ even}) \text{ or } (i_p < j_p \text{ and } t \text{ odd}) \quad (42)$$

The objective of the lunar logistics MILP is to minimise the total mass of goods launched into a translunar trajectory (arc [0,1]). The constraint explanations are as follows:

1. A: The net change in vehicles (and their associated stacks) entering and leaving a node are limited by supply (positive d), or demand (negative d).
   B: Other commodity types have fixed (f) and mutable (m) parts. Mutable demands are mapped onto the timeline of the fixed demands by the binary scheduling variable $B_{p,t}$. Goods can freely move between vehicles.
2. A: Sum of commodity mass carried by the vehicles travelling along an arc are limited by the payload capacities of the vehicles.
   B: ISRU infrastructure does not contribute towards vehicle payload capacity on the lunar surface, as it is deployed here.
   C & D: Vehicle oxidiser and fuel capacities, respectively.
3. A: Crew supplies consumption.
   B: ISRU infrastructure maintenance supply consumption on the lunar surface.
   C: Maintenance supplies are conserved in other locations.
   D: Oxidiser consumption across arcs is affected by trajectory mass fraction Z and boil-off β.
   E: Oxidiser can be refuelled by ISRU infrastructure on the lunar surface.
   F: Fuel consumption.
   G: Other commodities are simply conserved.
4. Commodities can only flow along allowed arcs, determined by binary parameter ε.
5. A: Each payload must be assigned a single launch time.
   B: Payloads must launch at the same time as their designated co-manifested payloads $\mathcal{C}_p$.
   C: Payloads must launch after or with their designated "soft" pre-cursor payloads $\mathcal{P}_p$.
   D: Payloads must launch after their designated "strict" pre-cursor payloads $\mathcal{Q}_p$.
   E: Payloads can only launch within their allowed launch windows, and on time indices aligning with their expected direction of travel.

Table 7: Parameter values used in the lunar logistics case study.

| Symbol | Parameter | Value |
|---|---|---|
| $\beta^f$ | Fuel boil-off (CH4) | 0.08% / day |
|  | Fuel boil-off (H2) | 0.1% / day [31] |
|  | Fuel boil-off (Storable) | 0 |
| $\beta^{ox}$ | Oxidiser boil-off (O2) | 0.025% / day [31] |
|  | Oxidiser boil-off (Storable) | 0 |
| $\phi$ | Propellant mixture ratio (CH4/O2) | 3.6 [32] |
|  | Propellant mixture ratio (H2/O2) | 6 |
|  | Propellant mixture ratio (Storable) | 2.61 [33] |
| $k$ | Crew Supply Consumption Rate | 8.7 kg/crew/day |
| $\mu$ | ISRU infrastructure maintenance supply consumption rate | 10 % infrastructure mass / year |
| $\rho$ | ISRU O2 production rate | 1.53 g propellant / kg infrastructure / day [34] |

The vehicles and payloads that define the parameters and requirements of the extended Artemis lunar surface exploration campaign are detailed in Tables 5 and 6 respectively, based on data from [32], [35], [36], [37], [38], [39], [40]. The remaining parameter values that govern the commodity dynamics in the model are listed in Table 7.

The lunar logistics problem is optimised using the time indices of the binary scheduling variables as the generated variable set, and the time index as the grouping index. Equation 6 is used to produce an initial feasible set of variables, with 18 time index groups selected. Due to the structure of the partially-static time-expanded network, it is necessary to pair time index groups together: an even (outbound direction) time step must be generated with the following odd (inbound direction) time step. Therefore, the group prices of the even time indices and their following odd indices are summed together.

## 5. Results and Discussion
### 5.1 Case Study 1

After generating 8 additional variables, the MILP is evaluated. The best-found solution to the Mars crewed mission ConOps have an objective (total propellant and tank mass) objective of 69746 kg, which closes the gap to the best bound found for the full problem. The details of the result are summarized in Table 9.

In summary, the spacecraft components are launched across 6 New Glenn launches. Note that in this case study, cost of launch is not considered, so the solver only seeks out feasible launch vehicles without considering their cost.

A 2000 km altitude parking orbit was selected as the assembly point for the spacecraft.

4 out of the maximum of 12 tanks are utilised, with 2 of them dropped after the initial departure from Earth to the trans-Mars trajectory. One more tank is dropped after the Mars capture burn. The trans-Earth injection and Earth capture manoeuvres are both carried out using propellant from the final remaining tank.

### 5.2 Case Study 2

The restricted MILP of Case Study 2 was first evaluated with just the minimum feasible set of variables. It was then re-evaluated after every 5 column generation algorithm iterations. The prices of the unused variables after the first iteration of the algorithm are shown in Figure 7. After adding the minimum variable set, and performing 5 column generation iterations, the solution improves to 576259 kg.

This solution to the lunar logistics problem is summarised in Table 10 It can be seen that the most mass-efficient of the explored solutions is one that maximises ride-sharing of payloads, with no dedicated launches for supporting payloads at all. Intuitively, this makes sense because it minimizes the amount of logistics vehicle and propellant that must be launched to deliver said payloads. Additionally, the sample return payloads are also combined together, with some Artemis crews bringing back several samples, and others bringing back none at all.

Further indices are added to the problem, but at this point the solve time begins to dramatically increase with no significant improvement to the solution.

### 4.3 Column Generation Algorithm Performance

Overall, the column generation algorithm performs strongly for the case studies detailed here. Figure 6 shows, for Case Study 1, the evolution of the best-found objective versus compute time for restricted problems of varying sizes compared to the full model. In all cases, the restricted problems find better solutions in with less computing time than the full model did.

Table 8 compares the best-found objective of Case Study 2 for different solution methods, and the amount of time taken to find said solution. Comparison is also made to the application of a previously-developed metaheuristic method [15].

### 3.4 Notes on the Choice of the Generated Variable Set

As discussed in Secs. 3.1 and 3.2, a key idea behind the proposed column generation method is to split the variables into a set that is always included and a set that is generated as needed. In this section, we reflect on the case studies, discuss how this split impacts the algorithm's performance, and retrospectively assess our particular choice of the generated variable sets and grouping indices.

The choice of the generated variable set and grouping index is driven by the pricing algorithm. In the case studies presented above, the pricing program involved calculating the reduced costs of continuous relaxations of the integer variables. As discussed earlier, this method of pricing is heuristic in nature. It is therefore important to consider the physical meaning of relaxing the integer variables, if they have any at all.

Table 8: Comparison of best-found objectives and the time taken to find them for different solve methods.

| Method | Objective (kg) | Compute Time (s) |
|---|---|---|
| Full MILP | 590992 | 238730 |
| Restricted problem: Minimum feasible $\gamma$ | 578305 | 6417 |
| Restricted problem: Min. $\gamma$ + 5 generated indices | 576259 | 87441 |
| Metaheuristics | 741457 | *not measured* |

Table 9: Results of the Mars crewed mission ConOps optimisation.

| Payload | Launch Vehicle | Tank Design Sample Capacity and Dry Mass (kg) | Propellant Usage for Manoeuvre # (kg) | | | |
|---|---|---|---|---|---|---|
| Assembly Orbit Altitude: | | 2000 km | 0 | 1 | 2 | 3 |
| Tank 0 | New Glenn Launch 1 | Capacity = 23100 Dry Mass = 2720 | 21630 | - | - | - |
| ~~Tank 1~~ | - | - | - | - | - | - |
| Tank 2 | New Glenn Launch 2 | Capacity = 5000 Dry Mass = 780 | 0 | 9830 | - | - |
| ~~Tank 3~~ | - | - | - | - | - | - |
| ~~Tank 4~~ | - | - | - | - | - | - |
| ~~Tank 5~~ | - | - | - | - | - | - |
| ~~Tank 6~~ | - | - | - | - | - | - |
| ~~Tank 7~~ | - | - | - | - | - | - |
| Tank 8 | New Glenn Launch 3 | Capacity = 19490 Dry Mass = 2370 | 18380 | - | - | - |
| ~~Tank 9~~ | - | - | - | - | - | - |
| ~~Tank 10~~ | - | - | - | - | - | - |
| Tank 11 | New Glenn Launch 4 | Capacity = 18280 Dry Mass = 2250 | 0 | 1020 | 7840 | 3340 |
| Spacecraft structure | New Glenn Launch 5 | - | | | | |
| Inert Payload | New Glenn Launch 6 | - | | | | |

Table 10: Launch time indices of the crewed missions of the extended Artemis program, and the launch time of the supporting payloads.

| Artemis Mission # | | Artemis Launch Time Index | Co-Manifested Payloads | Pre-launched Support Payloads |
|---|---|---|---|---|
| 7 | Outbound | 0 | Power Plant | - |
|   | Return | 1 | - | |
| 8 | Outbound | 10 | Habitat | - |
|   | Return | 12 | - | |
| 9 | Outbound | 17 | - | - |
|   | Return | 24 | Sample Return x2 | |
| 10 | Outbound | 36 | Pressurised Rover x2 | - |
|    | Return | 40 | Sample Return x2 | |
| 11 | Outbound | 48 | Fission Power Plant, Habitat | - |
|    | Return | 54 | - | |
| 12 | Outbound | 59 | - | - |
|    | Return | 66 | - | |
| 13 | Outbound | 69 | - | - |
|    | Return | 72 | Sample Return x4 | |
| 14 | Outbound | 84 | - | - |
|    | Return | 86 | - | |

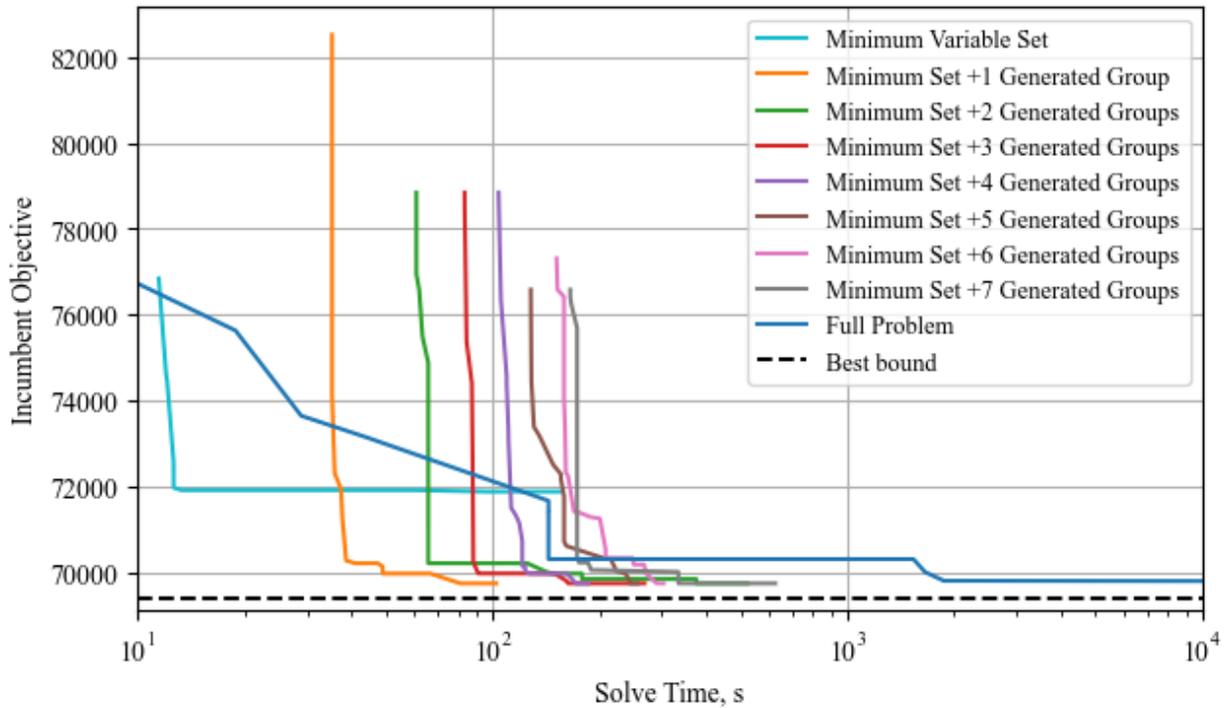

Figure 6: Best-found objective versus compute time of the crewed Mars mission mixed-integer model, for varying numbers of generated $\ell_{v,o}$ variables, ranging from the minimum feasible set of variables to the minimum set + 8 generated variable groups.

To assess our choice of the generated variable set, we evaluate several alternative options for the generated variable set in Case Study 1. Figures 8-10 compare the pricing of some different options, namely specifying the generated variable set $\Gamma$ to be the tank decision variable samples $x_{t,n,s,v}$ with grouping index $s$; $x_{t,n,s,v}$ with grouping index $v$; or the launch vehicle selection variable in the vehicle index $\ell_{v,o}$ with grouping index $v$. Note that the last one is the one chosen for the case study. As can be seen, not all of these options are good selections for the generated variable set. It can be seen in Figure 8 that the initially selected design was the extreme largest design. It happens that this is the most mass-efficient design. Relaxing $x$ effectively allows for fractional tanks. Therefore, the clear solution to the relaxed problem is to use fractional tanks scaled to the required size with the most efficient design. This of course isn't allowed in the integer problem: the tank fraction is required to match the design of the tank size. Therefore, the prices of the relaxed problem will not capture the value of smaller, but less mass-efficient, tank designs. In Figure 9, all of the prices are the same. This is because, in this specific formulation, there are no coefficients on $x$ or in its corresponding $b$ or $c$ vector entries that involve any launch vehicle-related parameters. Those parameters only act as coefficients on $\ell_{v,o}$, the prices of which are shown in Figure 10. Therefore, the prices of $x$ will not capture the usefulness of different launch vehicle types. Thus, our decision of using $\ell_{v,o}$ with grouping index $v$ is shown to be the most effective for this particular problem.

The choice of the generated variable set also influences the algorithm speed. The dimensions of $\ell_{v,o}$ are much smaller than $x_{t,n,s,v}$, so there are far fewer prices to compute per iteration of the process laid out in Figure 1.

5. Conclusions

The paper proposes and demonstrates that column generation can be a useful method for finding strongly performing solutions to complex exploration mission design problems, including problems that feature large numbers of categorical decisions. Provided that the set of generated variables is selected carefully, with attention paid to the structure of the problem and the meaning of the relaxation of integer variables, the restricted problem can produce solutions of equal or greater quality to the full problem in a much shorter amount of time.

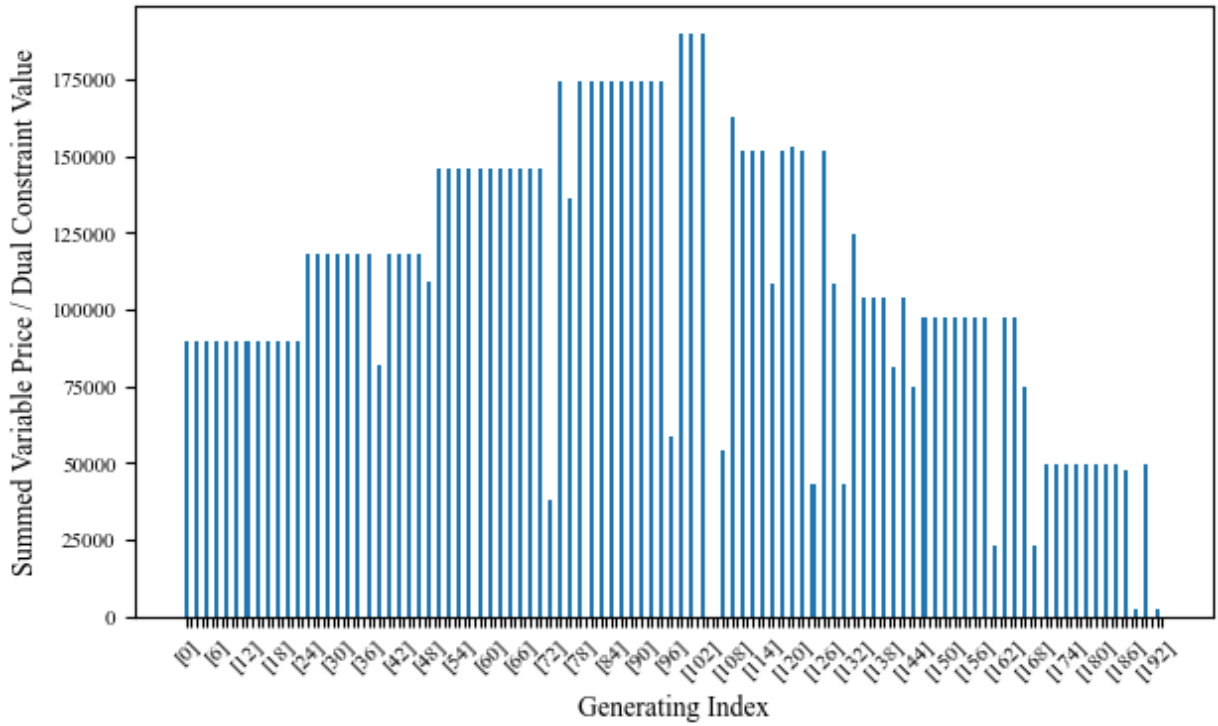

Figure 7: Total generating variable $B_{p,t}$ prices for each pair of grouping indices $t$ and $t+1$ $\forall$ *even* $t$ in Case Study 1.

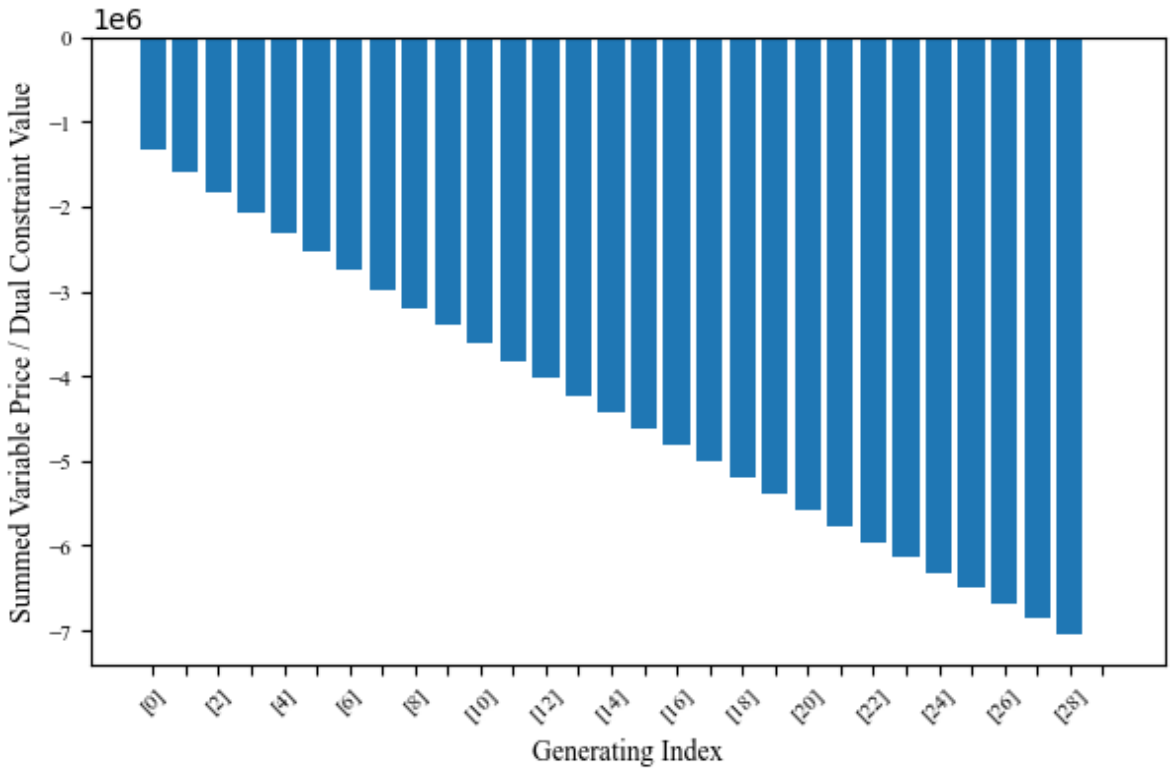

Figure 8: Prices of the generated variable set $x_{t,n,s,v}$ with grouping index $s$ after one column generation iteration in Case Study 2.

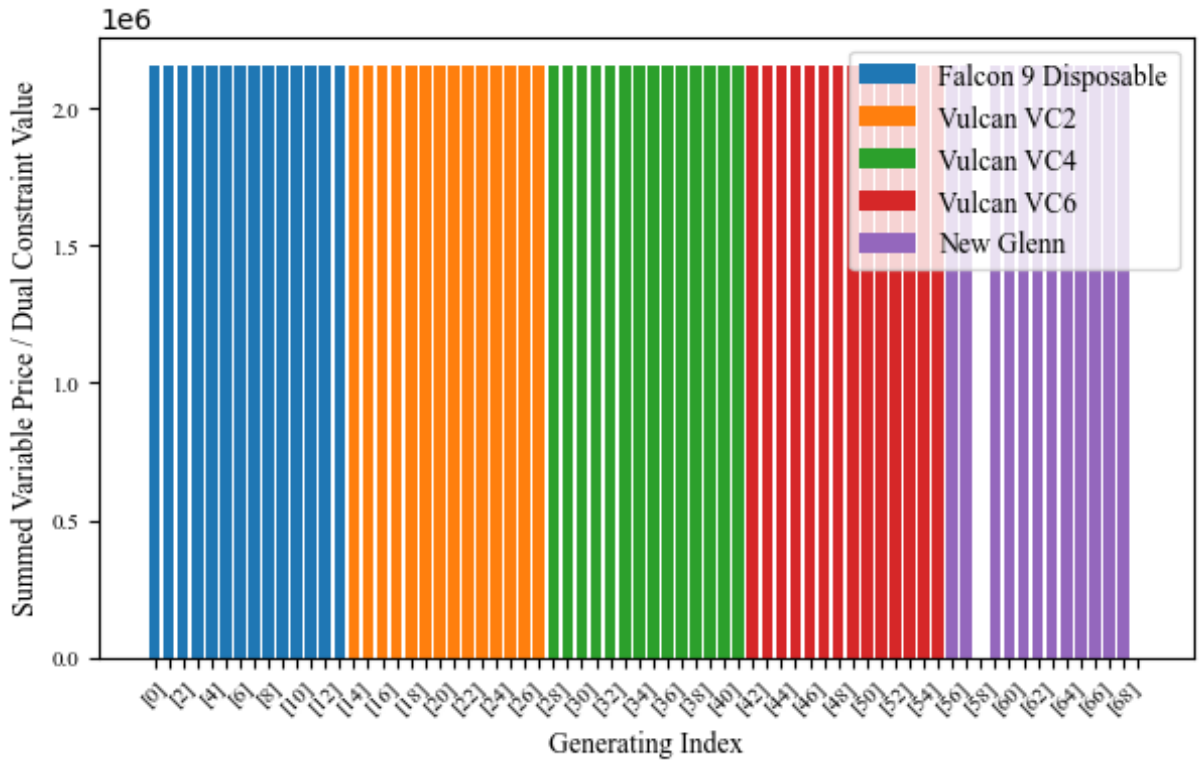

Figure 9: Prices of the generated variable set $x_{t,n,s,v}$ with grouping index $v$ after one column generation iteration in Case Study 2.

Figure 10: Prices of the generated variable set $\ell_{v,o}$ with grouping index $v$ after one column generation iteration in Case Study 2.

The two case studies discussed in this paper are formerly solved by fully sampling the discrete decision space, or by using metaheuristics. The application of mixed-integer programming with column generation to these problems reduces the computational time and resources required to address these design problems.

The algorithm presented here used the Dantzig-Wolfe decomposition to allow application to a wide range of ConOps optimisation problems. A useful focus of future work would be to formulate problem-specific pricing problems to more quickly converge to good solutions.

**Acknowledgment**

This research was supported by Jacobs Engineering Services and Science Capability Augmentation (ESSCA) contract through the Marshall Space Flight Center Advanced Concepts Office. We would like to express our gratitude to Dr. Stephen Edwards, Dr. Manuel Diaz, and Dr. Stephanie Zhu at the Advanced Concepts Office of NASA Marshal Space Flight Center for their invaluable feedback and assistance in this research.